\documentclass[titlepage,12pt]{article}
 \usepackage{amsmath,amsfonts,amsthm,amscd,amssymb,euscript}
 \usepackage{epsfig}
\def\Z{{\mathbb Z}}
\def\H{{\mathbb H}}

\newcommand{\calO}{{\mathcal{O}}}
\newtheorem{Th}{Theorem}[section] 
\newtheorem{Def}{Definition}[section]
\newtheorem{Lm}{Lemma}[section] 
\newcommand{\Prob}{{\hbox{Prob}}}

\newcommand{\genus}{{\hbox{\rm genus}}}
\newcommand{\diam}{{\hbox{\rm diam}}}

\begin{document}
\title{On the Genus of a Random Riemann Surface}
\author{Alexander
 Gamburd \\  MSRI\\
1000 Centennial Drive\\ Berkeley, CA 94720\\
gamburd@msri.org\\
\\
\\ Eran Makover \thanks{Partially supported by the NSF Grant DMS-0072534} \\ Department of Mathematics\\
Dartmouth College\\ Hanover, New Hampshire 03755\\ Eran.Makover@dartmouth.edu}
\date{December 2000}
\maketitle
\begin{abstract}
Brooks and Makover introduced an approach to random Riemann 
surfaces based on associating a dense set of them -- Belyi
surfaces -- with random cubic graphs.  In this paper, using 
Bollobas model for random regular graphs, we examine the
topological structure of these surfaces, obtaining in particular
an estimate for the expected value of their genus.
\end{abstract}

\section{Introduction}\label{sec.intro}

 In \cite{BM} Brooks and Makover constructed Riemann surfaces from
 oriented cubic graphs.  For each orientated graph $(\Gamma, \calO)$ they
 associate two Riemann surfaces, $S^O(\Gamma, \calO)$ a finite area
 noncompact surface, and $S^C(\Gamma, \calO)$ a compact surface.  The
 surface $S^O(\Gamma, \calO)$ is an orbifold cover of $\H/PSL(2,\Z)$
 and therefore shares some of the global geometric properties with the
 graph $(\Gamma, \calO)$.  The compact surface $S^C(\Gamma, \calO)$ is
 a conformal compactification of $S^O(\Gamma,\calO)$; Brooks and
 Makover proved that almost always the global geometry of $S^C(\Gamma,
 \calO)$ is controlled by the geometry of $S^O(\Gamma,\calO)$. 
 Moreover, according to Belyi theorem \cite{Bel} the surfaces $S^C(\Gamma,
 \calO)$  are precisely the Riemann surfaces
 which can  be defined over some number field and 
form a dense set in the space of all Riemann surfaces. 
  Therefore these surfaces, called Belyi surfaces, can be used
 to model a process of picking random Riemann surfaces, by picking a
 random graph with a random orientation.

In their work, Brooks and Makover use the graphs to get information on
the global geometry of the surfaces $S^{C}(\Gamma,\calO)$, but they do
not have control on the distribution of the surfaces in the
Teichm\"{u}ller spaces.  In this paper we investigate the topology of
Belyi surfaces, which is the first step in understanding the
distribution of $S^{C}(\Gamma,\calO)$ in different Teichm\"{u}ller
spaces.  Our main result (with $n$ denoting the number of vertices
in a cubic graph) is the following theorem:

\begin{Th}\label{thm:main}
 $$\boldsymbol{E}(\genus( S^C(\Gamma_{n}, \calO))) \sim O(n)$$
\end{Th}

The genera of these surfaces can be calculated by the Euler
formula where the number of vertices and edges is determined by the
size of the graph $\Gamma$.  The faces, which correspond to the cusps
of $S^O(\Gamma,\calO)$ can be described in a purely combinatorial way
from the oriented graph $(\Gamma,\calO)$.  Using Bollobas model for
random regular graphs \cite{Bol1} \cite{Bol2}, we estimate the number
of faces:
\begin{Th}\label{thm:face}
   There exist constants $C_1$ and $C_2$ such that, for a  large enough
   $n$:

 $$C_{1}\log(n)\le \boldsymbol{E}(\mathcal{F}(\Gamma_{n},\calO)) \le
 C_{2} \sqrt{n}$$
 Where $\mathcal{F}$ is the number of faces of the oriented graph.

\end{Th}

 The Euler formula then gives our main result.

We will start in section \ref{sec.Bey} in a review on the construction
of the Belyi surfaces from cubic graphs, and we will show how to
calculate the genus of the surfaces.  In section \ref{sec.CE} we will
review the Bollobas model and use it to calculate the expected value
of the number of faces in Belyi surfaces.

Note: Our main interest is the construction of the Belyi surfaces,
this is why we work with cubic graphs.  Other authors \cite{sh} work
on the average value of the genus of a given graph: our result can be
easily generalized to $k$-regular graphs.

\begin{figure}[!h]
\begin{center}

\epsfig{file=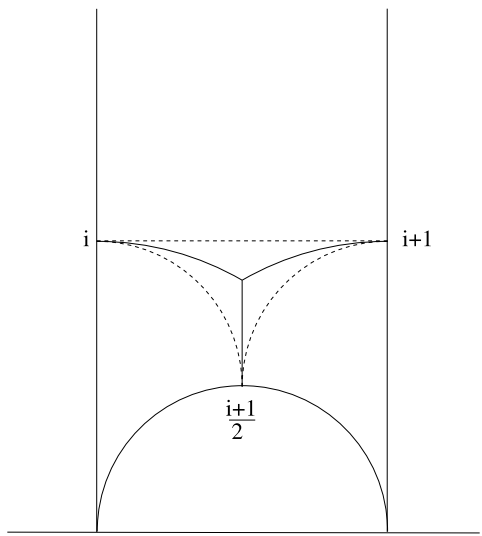, height=7cm}
\caption{The marked ideal triangle $T$}\label{triangle}
\end{center}
\end{figure}
\section{Belyi surfaces}\label{sec.Bey}
An orientation $\calO$ on a graph $\Gamma$ is a function which
associates to each vertex $v$ a cyclic ordering of the edges emanating
from $v$.  We build the surface $S^O(\Gamma,\calO)$ from an oriented
graph as follows: we take the ideal hyperbolic triangle $T$ with
vertices $0, 1$, and $\infty$ shown in Figure \ref{triangle}.  The
solid lines in Figure \ref{triangle} are geodesics joining the points
$i, i+1$, and ${\frac{i+1}{2}}$ with the point ${\frac{1 + i
\sqrt{3}}{2}}$, while the dotted lines are horocycles joining pairs of
points from the set $\{i, i+1, {\frac{i+1}{2}}\}$.  We may think of
these points as ``midpoints'' of the corresponding sides of the ideal
triangles, even though the sides are of infinite length.  We may also
think of the three solid lines as edges of a graph emanating from a
vertex.  We may then give them the cyclic ordering $(i, i+1,
{\frac{i+1}{2}})$.  Given a cubic graph with orientation
$(\Gamma,\calO)$ we build a non-compact Riemann surface denoted by
$S^O(\Gamma, \calO)$, by associating to each vertex an ideal triangle,
and gluing neighboring triangles.  We glue two copies of $T$ along the
corresponding sides, subject to the following two conditions:

\begin{description}

\item[a] the ``midpoints'' of the two sides are glued together,

and

\item[b] the gluing preserves the orientation of the two copies of
$T$.

\end{description}

The conditions (a) and (b) determine the gluing uniquely.  It is
easily seen that the surface $S^O(\Gamma, \calO)$ is a complete
Riemann surface with a finite area equal to $\pi n$, where $n$ is the
number of vertices of $\Gamma$, and that the horocycles of the copies
of $T$ fit together to give closed horocycles about the cusps of
$S^O(\Gamma, \calO)$.

We denote by $S^C(\Gamma, \calO)$  the conformal compactification of $S^O(\Gamma,
\calO)$. Using the results about random cubic graphs
Brooks and Makover proved \cite{BM}:
\begin{Th}\label{typical}
There exist constants $C_1$, $C_2$, $C_3$, and $C_4$ such that, as $n
\to \infty$:

\begin{description}
 \item[a] The first eigenvalue $\lambda_1(S^C(\Gamma, \calO))$
 satisfies $$\Prob_n[ \lambda_1(S^C(\Gamma, \calO)) \ge C_1] \to 1.$$
 
 \item[b]The Cheeger constant $h(S^C(\Gamma, \calO))$ satisfies
 $$\Prob_n [h(S^C(\Gamma,\calO)) \ge C_2] \to 1.$$

\item[c] The shortest geodesic $\rm{syst}(S^C(\Gamma, \calO))$
satisfies $$Prob_n[\rm{syst}(S^C(\Gamma, \calO)) \ge C_3] \to 1.$$

\item[d] The diameter $\diam(S^C(\Gamma, \calO))$ satisfies
$$\Prob_n[\diam(S^C(\Gamma, \calO)) \le C_4 \log(\genus(S^C(\Gamma,
\calO)))] \to 1.$$
\end{description}
\end{Th}

The topology of the surface can be read from the oriented graph using
left-hand-turn paths.  A {\emph left-hand-turn path} on $\Gamma$ is a closed path
on $\Gamma$ such that, at each vertex, the path turns left in the
orientation $\calO$.

On $S^{O}(\Gamma,\calO)$ a left-hand path describe a closed path
around a cusp.

Let $n:=\mathcal{V}(\Gamma)$ and $l(\Gamma,\calO)$ be the number of
disjoint left-hand paths.  For a cubic graph the  number of edges is
$\mathcal{E}(\Gamma)=\frac{3n}{2}$ and the number of faces is
$\mathcal{F}=l(\Gamma,\calO)$. Therefore we can write the Euler
formula:

\begin{equation}
\label{e:gen}
genus(S^{C}(\Gamma_{n},\calO))=genus(S^{O}(\Gamma_{n},\calO))=1+\frac{n-l}{2}.
\end{equation}

\begin{figure}[!h]
\begin{center}

\epsfig{file=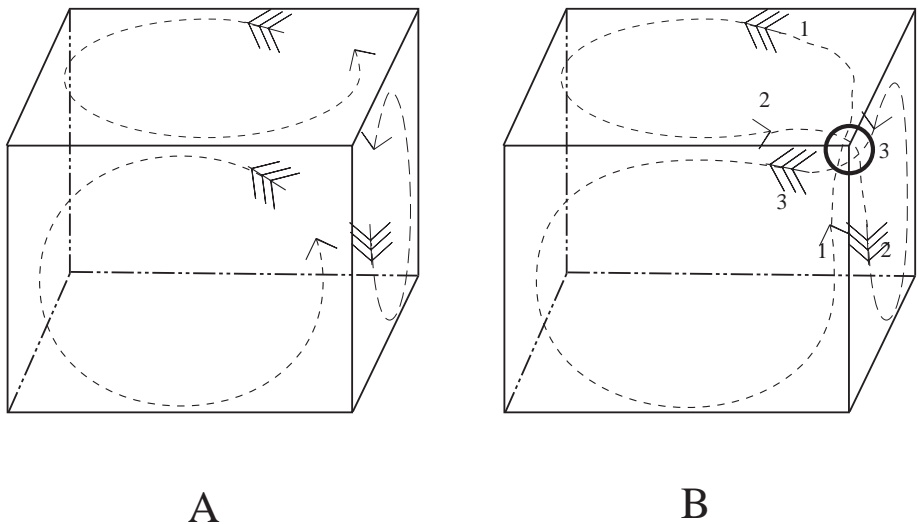, height=7cm}
\caption{Changing orientation on the cube }\label{cube}
\end{center}
\end{figure}

To estimate the genus of a random surface we need to estimate the
number of faces (or left-hand paths of the oriented graph).  It is
important to observe that a left-hand path is not necessarily a simple
closed path on the graph.  For example, if we take the 1-skeleton of
the cube, with the usual orientation
(Figure \ref{cube} \textbf{A}),
we have $\mathcal{F}=6$, and all the faces are simple paths of the
graph.
In Figure \ref{cube} \textbf{B} we changed the orientation of
the right upper vertex.  With the new orientation
Now if we change the orientation of the right upper vertex,
 the three simple
faces that were adjacent to the upper right vertex before we changed
orientation are now joined to one composite face while the other three
faces are unchanged hence $\mathcal{F}=4$.

In section \ref{sec.CE} we obtain an upper and lower bound 
on the number of faces given in theorem \ref{thm:face}. The lower bound 
is obtained by counting the number
of faces that are simple closed paths on the graph.

The upper bound on the number of faces is obtained by  first dividing the faces
into two groups:

\begin{description}
    \item[a] Small faces $\hbox{\rm length}(\mathcal{F})\le \sqrt(n)$

    \item[b] Large faces $\sqrt{n} \le \hbox{\rm length}(\mathcal{F})$
\end{description}

The number of large faces is easily bounded: since the total number of edges
is $\frac{3n}{2}$ and each edge has two sides,  we obtain that 
$\#(\mathcal{F}_{\ge\sqrt{n}})\le 3 \sqrt{n}$.

To estimate the number of small faces we first introduce the notion 
of root.

\begin{Def}\label{def:root}
   A root in $(\Gamma,\calO)$ is a simple closed path in $\Gamma$,
   such that the orientation $\calO$ agrees with the root in all but
   maybe one vertex.
\end{Def}
\begin{Lm}
Each face contains at least one root: 
\end{Lm}
\noindent
{\bf{Proof}}:\quad
Pick a vertex $v_{1}$ and start
walking along a left hand path until the path intersects itself, i.e.,
$v_{i}=v_{j}$ $i\ne j$ is the first intersection, then the cycle
$v_{i}, \cdots v_{j}$ is a root (in all the vertices $v_{l}$ $i<l<j$
the orientation agrees with the cycle). \\

In section \ref{sec.CE} we will prove that the expected 
number of roots of length less 
than  $\sqrt{n}$, and therefore the total number of small faces, is bounded
by  $O (\sqrt{n})$.
By adding the number of short and  long faces we get that the
total number of faces is $< O(\sqrt{n})$; from equation
\ref{e:gen} we deduce our main result \ref{thm:main}.

\section{Computing the Expectation}\label{sec.CE}
This section is devoted to the proof of theorem \ref{thm:face}.
Let $G_{n, d}$ be the space of $d$-regular graphs with $n$ vertices.  We
shall, following Bollobas, represent our graphs as the images of
so-called configurations.  Let $W=\bigcup_{j=1}^{n} W_{j}$ be a fixed
set of $2 m = nd$ vertices, where $| W_{j}|=d$.  A
\emph{configuration} $F$ is a partition of $W$ into $m$ pairs of
vertices called \emph{edges} of $F$.  Clearly there are
\begin{equation} \label{e:a1}
N=N(m)=\binom{2m}{2} \binom{2m-2}{2} \dots \binom{2}{2} / m! 
=\frac{(2m)_{m}}{2^{m}}
\end{equation}
configurations.  (We write $(a)_{b}=a(a-1) \dots (a-b+1)$).

Let $\Phi$ be a set of configurations.  We now define a map $\Phi
\rightarrow G_{n, d}$ as follows.  Given a configuration $F$ let
$\phi(F)$ be the graph with vertex set $V={1, \dots, n}$ in which $ij$
is an edge iff $F$ has a pair with one end in $W_{i}$ and the other in
$W_{j}$.  Every $G \in G_{n, d}$ is the image of $G=\phi(F)$ for
$(d!)^{n}$ configurations.  The number of configurations containing a
given fixed set of $l$ edges is
\begin{equation} \label{e:a2}
\begin{split}
N_{l}(m) &= \binom{2m-2l}{2} \binom{2m-2l-2}{2} \dots \binom{2}{2}
/(m-l)!  \\
&=\frac{N(m)}{(2m-1)(2m-3) \dots (2 m -2l +1)} \\
\end{split}
\end{equation}

\vskip 3pt

From (\ref{e:a1}) and (\ref{e:a2}), the probability that a
configuration contains a given set of $l$ edges is
\begin{equation} \label{e:a3}
\begin{split}
\frac{N(m-l)}{N(m)} &= \frac{1}{(2m-1)(2m-3) \dots (2m-2l+1)} \\
&= \frac{1}{2^{l}} \frac{\Gamma (m-l+\frac{1}{2})}{\Gamma
(m+\frac{1}{2})}\\
\end{split}
\end{equation}

where we use
\[
1 \cdot 3 \cdot 5 \cdot \dots (2 n -1) = 2^{n} \Gamma (n+\frac{1}{2})
\cdot \frac {1}{\sqrt{\pi}}
\]

For $k \in \bf{N}$ a \emph{k-cycle} of a configuration is a set of $k$
edges, say ${e_1, \dots , e_{k}}$ such that for some $k$ distinct
groups $W_{j_{1}}, \dots , W_{j_{k}}$ the edge $e_{i}$ joins
$W_{j_{i}}$ to $W_{j_{i+1}}$, where $W_{j_{k+1}} \equiv W_{j_{1}}$.  A
1-cycle is said to be a \emph{loop} and a 2-cycle is a
\emph{coupling}.  Given a configuration $F$, denote by $X_{k}(F)$ the
number of k-cycles.  If we are to restrict to graphs in $G_{n, d}$ to
have no loops and no multiple edges, then not every $\phi(F)$ belongs
to $G_{n, d}$ but only those satisfying $X_{1}(F) =X_{2}(F) =0$; such
graphs are called simple.  (In order to restrict to graphs with
$X_{1}=X_{2}=0$ we have to multiply all the ensuing estimates by
$\bf{P}(\rm{simple}) \sim \exp(\frac{1-d^{2}}{4})$.)

Let $C_{k}$ be the number of sets of pairs of vertices in $V$ which
can be k-cycles of configurations.  By elementary counting:
\begin{equation} \label{e:a4}
\ C_{k} = \frac{1}{2 k}\frac{n!}{(n-k)!} (d (d-1))^{k}
\end{equation}

Using (\ref{e:a3}) and (\ref{e:a4}) we get the following expression
for the expected number of $k$-cycles:

\begin{equation} \label{e:a5}
\boldsymbol{E}(X_{k}) = \frac{1}{2k} \frac{n!}{(n-k)!} (d (d-1))^{k}
\frac{1}{2^{k}} \frac{\Gamma (m-k+\frac{1}{2})}{\Gamma
(m+\frac{1}{2})}
\end{equation}

To obtain a lower bound on the number of faces $\mathcal{F}$, we count
the number of faces that are simple closed paths.  For each simple
closed path of length $k$ the probability of correct orientation is
$\frac{1}{(d-1)^{k}}$, consequently we have:

\begin{equation} \label{e:a9}
\boldsymbol{E}(\mathcal{F}) > \sum_{3}^{m} \frac{\boldsymbol{E}(X_k)}
{(d-1)^k} =\sum_{3}^{m}\frac{a_k}{k},
\end{equation}

where
\begin{equation}   \label{e:a7}
a_{k} = \frac{d^{k}}{2^{k}} \frac{n!}{(n-k)!} \frac{\Gamma
(m-k+\frac{1}{2})}{\Gamma (m+\frac{1}{2})}.
\end{equation}

Now using Stirling's formula (with $x$ positive real),
\[
\Gamma(x) =\sqrt{2 \pi}x^{x-\frac{1}{2}}e^{-x}\left(1+
O\left(\frac{1}{x}\right)\right),
\]

we have for $k=o(n^{2/3})$
\[
\frac{n!}{(n-k)!} = 
n^{k} e^{\frac{k^{2}}{2n}} \left(1+O\left(\frac{k}{n}\right) +
O\left(\frac{k^{3}}{n^{2}}\right) \right);
\]

and similarly
\[
\frac{\Gamma (m-k+\frac{1}{2})}{\Gamma (m+\frac{1}{2})}= m^{-k}
e^{\frac{k^{2}}{2m}} \left(1+O\left(\frac{k}{n}\right)\right).
\]

Recalling that $m=\frac{nd}{2}$, we have
\begin{equation}  \label{e:a8}
a_{k} = e^{-\frac{(d-2)k^{2}}{2 n d}}\left(1+O\left(\frac{k}{n}\right)\right).
\end{equation}

The estimate (\ref{e:a8}) and the fact that $a_k$ are monotonic
decreasing in $k$, show that if we let $k$ range from $3$ to $m^t$
where $t>1/2$, the remaining terms in (\ref{e:a9}) are negligible.
With that in mind, let $b$ be an arbitrary positive number 
independent of $m$.
If $k^2 \leq bm$, the estimate in (\ref{e:a8}) is between 
$1$ and $\exp{(-b/2)}$.  Using 
\[
1+\frac{1}{2} + \dots + \frac{1}{l} \sim \log{l}
\]
we get that the asymptotic estimate for that part of sum 
in (\ref{e:a9}) with $3 \leq k \leq \sqrt{bm}$ is between 
$(\log m)/2$ and $e^{-b/2} (\log m)/2$.  For the remainder 
of the sum, $\frac{1}{k} < \sqrt{bm} $, and the sum of 
$\exp(-k^2/2m)$ is asymptotic to $\sqrt{m}$ times an integral of 
$\exp(-x^2/2)$.  It follows that for any value of 
$b$, the contribution of this part of sum is bounded as $m \to \infty$,
and so (\ref{e:a9}) is asymptotic to a function 
of $m$ which lies between $(\log m)/2$ and $\exp(-b/2) (\log m)/2$.
Since $b$ is arbitrary, $\exp(-b/2)$ can be made as close to 1 as desired;
hence we obtain that the sum in 
 (\ref{e:a9}) is asymptotic to $(\log m)/2$.

\vskip 5pt 

Now we turn to estimating  the number of faces from above. Recall that the number
of large faces is bounded from above by $\sqrt{m}$, so we need to estimate from above
 the number of small faces.  As their length is less than $\sqrt{m}$, so is the length 
of their roots; thus the estimate on the number of roots of length less than $\sqrt{m}$
will suffice.

Now given a cycle of length $k$, we can turn it into a root in $(k+1)$ ways: one simple face 
and $k$ proper roots.  Furthermore, the probability of correct orientation is 
$\frac{1}{(d-1)^{k}}$. So we have the following estimate
(with $a_k$ given in \ref{e:a7}):

\[
\boldsymbol{E}(\mathcal{F}_{\le \sqrt{m}}) < \sum_{k=3}^{\sqrt{m}}
\frac{k+1}{(d-1)^k} \boldsymbol{E}(X_{k}) = \frac{1}{2}
\sum_{k=3}^{\sqrt{m}} a_{k} + 
\frac{1}{2}\sum_{k=3}^{\sqrt{m}} \frac{a_{k}}{k}.
\]

The second sum is clearly dominated by the sum in (\ref{e:a9}), which we 
established to be of order $O(\log{m})$.  For the first sum,
 using (\ref{e:a8}),
we obtain

\[
 \sum_{k=3}^{\sqrt{m}} a_{k} = \sum_{k=3}^{\sqrt{m}}
e^{-\frac{(d-2)k^{2}}{2 n d}}\left(1+O\left(\frac{k}{n}\right)\right) <
\sqrt{nd}\int_{0}^{\infty}e^{-\frac{x^{2}}{2}} \mathrm{d} x + O(1) =
O(\sqrt{n})
\]

We have therefore established upper and lower bound of 
$O(\sqrt{n})$ and $O(\log{n})$ respectively
for the number of faces, proving the 
theorem \ref{thm:face}.

\ifx\undefined\bysame
\newcommand{\bysame}{\leavevmode\hbox to3em{\hrulefill}\,}
\fi

\end{document}